\documentclass[12pt,leqno]{amsart}
\usepackage{amsmath,amssymb}
\usepackage[latin1]{inputenc}

\newtheorem{thm}{Theorem}
\newtheorem{lem}[thm]{Lemma}

\newcommand{\PP}{\mathcal{P}}
\newcommand{\LL}{\mathcal{L}}
\newcommand{\MM}{\mathcal{M}}
\newcommand{\HH}{\mathcal{H}}
\newcommand{\Ss}{\mathcal{S}}
\newcommand{\WW}{\mathcal{W}}

\newcommand{\lcm}{\operatorname{lcm}}

\newcommand{\ord}{\operatorname{ord}}
\newcommand{\oo}{\operatorname{ord^*}}

\newcommand{\bigA}{T}

\newcommand{\rmkd}[1]{}

\title[Period of the linear congruential and power generators]
{On the period of the linear congruential and power generators}

\author{P\"ar Kurlberg}
\email{kurlberg@math.chalmers.se}
\address{Department of Mathematics\\ 
Chalmers University of Technology\\
SE-412 96 Gothenburg  \\
Sweden}

\author{Carl Pomerance}
\email{carlp@math.dartmouth.edu}
\address{Mathematics Department\\
Dartmouth College\\
Hanover, NH 03755-3551\\
U.S.A.}

\thanks{P.K.\ supported in part by the
National Science Foundation (DMS 0071503), the Royal Swedish
Academy of Sciences and the Swedish Research Council.
C.P.\ supported in part by the National Science Foundation.}

\begin{document}



\maketitle

\section{Introduction}

We consider two standard pseudorandom number generators from number theory:
the linear congruential generator and the power generator.  For the
former, we are given integers $e,b,n$ (with $e,n>1$) and a seed $u=u_0$, and we compute the
sequence
\[
u_{i+1}~=~eu_i+b~(\text{mod }n).
\]
This sequence was first considered as a pseudorandom number generator by D.~H.~Lehmer.
For the power generator we are given
integers $e,n>1$ and a seed $u=u_0>1$, and we compute the sequence
\[
u_{i+1}~=~u_i^e~(\text{mod }n)
\]
so that $u_i=u^{e^i}~(\text{mod }n)$.  A popular case is $e=2$, which is called
the Blum--Blum--Shub (BBS) generator.

Both of these generators are periodic sequences, and it is of interest
to compute the periods.  To be useful, a pseudorandom number generator
should have a long period.  In this paper we consider the problem of
the period statistically as $n$ varies, either over all integers, or
over certain subsets of the integers that are used in practice, namely
the set of primes and the set of ``RSA moduli," that is, numbers which
are the product of two primes of the same magnitude.

If $(e,n)=1$, then the sequence $e^i~(\text{mod }n)$ is purely periodic and its
period is the least positive integer $k$ with $e^k\equiv1~(\text{mod }n)$.  We
denote this order as $\ord(e,n)$.  If $(e,n)>1$, the sequence $e^i~(\text{mod }n)$
is still (ultimately) periodic, with the period given by $\ord(e,n_{(e)})$ where
$n_{(e)}$ is the largest divisor of $n$ that is coprime to $e$.  (The aperiodic
lead-in to such a sequence has length bounded by the binary logarithm of $n$.)
In this paper we shall denote $\ord(e,n_{(e)})$ by $\oo(e,n)$.  The periods of both
the linear congruential and power generators may be described in terms of this
function.  For the linear congruential generator we have
$u_i=e^i(u+b(e-1)^{-1})-b(e-1)^{-1}~(\text{mod }n)$ when $e-1$ is
coprime to $n$, so that if we additionally have $u+b(e-1)^{-1}$ coprime to $n$,
the period is exactly $\oo(e,n)$.
In general, the period is always a divisor of $\oo(e,n)(e-1,n)$.  

For the power generator, the period is exactly $\oo(e,\oo(u,n))$. 
We shall assume that $u$ is chosen so that $\oo(u,n)$ is as large
as possible for a given modulus $n$.\footnote{At the end of the paper we briefly
consider the general case where this assumption is not made.}  This maximum is denoted
$\lambda(n)$, following Carmichael.  First described by Gauss, $\lambda(n)$
is the order of the largest cyclic subgroup of $({\bf Z}/n{\bf Z})^\times$.
It satisfies $\lambda([a,b])=[\lambda(a),\lambda(b)]$, where $[~,~]$
denotes the least common multiple.  Further, for a prime power $p^\alpha$
we have $\lambda(p^\alpha)=\phi(p^\alpha)=(p-1)p^{\alpha-1}$, except when
$p=2,\alpha\ge3$ in which case $\lambda(2^\alpha)=2^{\alpha-2}$.  
For the power generator, we thus will study $\oo(e,\lambda(n))$.  Note that it
is especially important to use the function $\oo$ rather than $\ord$ when
considering the modulus $\lambda(n)$, since for $n>2$, $\lambda(n)$ is
always even, and in general, $\lambda(n)$ is divisible by the fixed
number $e$ for a set of numbers $n$ of asymptotic density~1.

We begin by reviewing some of the literature on statistical properties of
$\oo(e,n)$.
In \cite{pappalardi-p-1-div} Pappalardi showed that there exist
$\alpha,\delta>0$ such that $\ord(e,p) \geq p^{1/2} \exp( (\log
p)^\delta)$ for all but $O(x/\log^{1+\alpha} x)$ primes $p \leq x$.
He also asserted, assuming the Generalized Riemann 
Hypothesis\footnote{More precisely, that the Riemann hypothesis holds for $L$-functions
associated with certain Kummer extensions} (GRH), that if
$\psi(x)$ is any increasing function 
tending to infinity as $x$ tends to infinity, then $\ord(e,p) > p/\psi(p)$ for
all but $O(\pi(x)\log(\psi(x))/\psi(\sqrt{x}) )$ primes $p \leq x$, where
as usual, $\pi(x)$ is the total number of all primes $p\le x$.
(Although stated for any unbounded
monotone function $\psi(x)$, it appears that the proof only supports
the case when $\psi(x)$ is increasing rather slowly.  A similar result
with $\psi(x) \le (\log x)^{1-\epsilon}$ is proved in the first author's
paper \cite{artincat}.  In Theorem~ \ref{thm:order-on-grh} we obtain a
small, yet for our purposes crucial,  strengthening of this result.)
In \cite{erdos-murty}, Erd\H os and Murty showed that
if $\epsilon(x)$ is any decreasing function tending to zero as $x$ tends
to infinity, then $\ord(e,p) \geq p^{1/2+\epsilon(p)}$
for all but $o(\pi(x))$ primes $p \leq x$, and in
\cite{indlekofer-timofeev-divisors-of-shifted-primes} Indlekofer and
Timofeev gave a similar lower bound with an explicit estimate on the
number of exceptional primes.  Further, it follows immediately
 from work of Goldfeld, Fouvry, and Baker--Harman that there is
a positive constant $\gamma$ such that
$\ord(e,p)>p^{1/2+\gamma}$ for a positive proportion of the primes $p$.

The period of the power generator $u^{e^i}~({\rm mod}~pl)$ was studied 
in Friedlander, Pomerance and Shparlinski \cite{pomerance-power-generator}, 
where $p,l$ are primes of the same
magnitude.  One of the results there is that this period is $>(pl)^{1-\epsilon}$
for most choices of $u,e,p,l$.  However, once the exponent $e$ is fixed,
say at 2, the results of \cite{pomerance-power-generator} are noticeably weaker.

As for $\ord(e,n)$ for $n$ a positive integer, in \cite{cat2} Kurlberg and
Rudnick proved that there exists $\delta>0$ such that $\ord(e,n) \gg
n^{1/2} \exp( (\log n)^\delta))$ for all but $o(x)$ integers $n \leq
x$ that are coprime to $e$.  Further, in \cite{artincat}, Kurlberg showed that the 
GRH implies that for each $\epsilon>0$, we have
$\ord(e,n) \gg n^{1-\epsilon}$ for all
but $o(x)$ integers $n\leq x$ that are coprime to $n$, and in 
\cite{li-pomerance-artin-for-composite} Li and Pomerance  improved
the lower bound to $\ord(e,n) \geq n (\log n)^{-(1+o(1))\log \log \log n}$,
a result that is best possible.

To complement these theorems we give some new results on
$\ord(e,n)$ and $\oo(e,n)$.
\begin{thm}
\label{thm:big-order-mod-n}
Results on $\oo(e,n)$:
\begin{enumerate}
\item
Suppose $\epsilon(x)$ tends to zero arbitrarily slowly as $x \to \infty$.
Then $\oo(e,n)\ge n^{1/2+\epsilon(n)}$
for all but $o_\epsilon(x)$ integers $n\le x$.
\item
There is a positive constant $\gamma_1$
such that $\ord(e,n)\ge n^{1/2+\gamma_1}$ for a positive proportion
of the integers $n$.
\end{enumerate}
\end{thm}
These relatively easy results, together with the GRH-conditional
results mentioned above, become the model for the principal results
of this paper.  We consider the power generator for 3 classes
of moduli: primes, the products of two primes of the same magnitude,
and general moduli.

\begin{thm}
\label{thm:big-order-mod-p-1}
Results on $\oo(e,p-1)$:
\begin{enumerate}
\item
Suppose $\epsilon(x)$ tends to zero arbitrarily slowly as $x \to \infty$.
Then $\oo(e,p-1) \geq p^{1/2+\epsilon(p)}$ for all but 
$ o_{\epsilon}(\pi(x)) $ primes $p \leq x$.  
\item
There is a positive constant
$\gamma_2$ such that $\oo(e,p-1)\ge p^{1/2+\gamma_2}$ for
a positive proportion of the primes $p$.  
\item
(GRH) For each fixed $\epsilon>0$ we have $\oo(e,p-1)>p^{1-\epsilon}$
for all but $o_\epsilon(\pi(x))$ primes $p\le x$.
\end{enumerate}
\end{thm}

Consider moduli $pl$ where $p,l$ are primes with $p,l\le Q$ (where $Q$ is
an arbitrary bound).
Using our results on $\oo(e,p-1)$, we can prove the following theorem.
\begin{thm}
\label{thm:main-epsilon}
Results on $\oo(e,\lambda(pl))$:
\begin{enumerate}
\item
Suppose $\epsilon(x)$ tends to zero arbitrarily slowly as $x \to \infty$.
Then $\oo(e,\lambda(pl)) \geq
(pl)^{1/2+\epsilon(pl)}$ for all but $o_{\epsilon}(\pi(Q)^2)$ pairs of
primes $p,l\le Q$.
\item
There is a positive constant
$\gamma_3$ such that for a positive proportion of the pairs of primes 
$p,l\le Q$, we have $\oo(e,\lambda(pl))\ge (pl)^{1/2+\gamma_3}$.
\item
(GRH) For each fixed $\epsilon>0$ we have $\oo(e,\lambda(pl))>(pl)^{1-\epsilon}$
for all but $o_\epsilon(\pi(Q)^2)$ pairs of primes $p,l\le Q$.
\end{enumerate}
\end{thm}

Instead of considering specifically RSA moduli $n=pl$,
one may consider the general case where no restriction is made on the modulus $n$.
As we have seen, the length of the period for the
sequence $(u_i)$ is bounded by $\oo(e,\lambda(n))$.  
In our last theorem we
establish similar results as above for this order.
\begin{thm}
\label{thm:lambda}
Results on $\oo(e,\lambda(n))$:
\begin{enumerate}
\item
Suppose $\epsilon(x)$ tends to zero arbitrarily slowly as $x \to \infty$.
Then $\oo(e,\lambda(n)) \geq
n^{1/2+\epsilon(n)}$ for all but $o_{\epsilon}(x)$ integers $n\le x$.
\item
There is a positive constant
$\gamma_4$ such that $\oo(e,\lambda(n))\ge n^{1/2+\gamma_4}$ for
a positive proportion of the integers $n$.  
\item
(GRH) For each fixed $\epsilon>0$ we have $\oo(e,\lambda(n))>n^{1-\epsilon}$
for all but $o_\epsilon(x)$ integers $n\le x$.
\end{enumerate}
\end{thm}

We actually achieve a best-possible result in part 3 of Theorem \ref{thm:lambda},
showing that, on assumption of the GRH, that
\[
\oo(e,n)~=~n\cdot\exp\left(-(1+o(1))(\log\log n)^2\log\log\log n\right)
\]
as $n\to\infty$ through a set of asymptotic density 1.

\medskip

\noindent{\bf Acknowledgement}.  We would like to thank
Igor Shparlinski for several helpful conversations.

\section{Preliminary ideas}
\label{sec:preliminaries}

In this section we present an argument that shows that
$\oo(e,n)>n^{1/2+\epsilon(n)}$ on a set of asymptotic density~1; that is,
we prove the first item in Theorem~\ref{thm:big-order-mod-n}.
This argument will then be a model for the analogous item in each of
Theorems~\ref{thm:big-order-mod-p-1}, \ref{thm:main-epsilon}, 
\ref{thm:lambda}.  

We begin with a useful lemma.  The proof appeared in \cite{cat2}, section 5.1,
but for completeness we give a somewhat shorter argument here.

\begin{lem}
\label{lem:ordineq}
For any natural number $n$ we have
\[
\oo(e,n)~\ge~\frac{\lambda(n)}{n}\prod_{p|n}\oo(e,p)
~=~\frac{\lambda(n)}{n}\prod_{p|n,~p\nmid e}\ord(e,p).
\]
\end{lem}
\begin{proof}
The equality is trivial.  For the inequality, note that for positive
integers $a_i,b_i$ 
we have
\[
\lcm\{a_1b_1,\dots,a_kb_k\}~|~b_1\cdots b_k\cdot\lcm\{a_1,\dots,a_k\},
\]
as each $a_ib_i$ divides $b_1\cdots b_k\cdot\lcm\{a_1,\dots,a_k\}$.
We apply this with the $a_i$'s being the various $\oo(e,p)$ for $p|n$ and
the corresponding $b_i$'s being $\lambda(p^\beta)/\oo(e,p)$, where
$p^\beta\|n$.  Then $\lcm\{a_1b_1,\dots,a_kb_k\}=\lambda(n)$.  Further,
$\oo(e,n)$ is divisible by $\lcm\{a_1,\dots,a_k\}$, so that
\[
\frac{\lambda(n)}{n}~\le~\frac{\oo(e,n)}{n}\prod_{p^\beta\|n}\frac{\lambda(p^\beta)}
{\oo(e,p)}~\le~\frac{\oo(e,n)}{\prod_{p|n}\oo(e,p)}.
\]
\end{proof}

Suppose $\mathcal{P}$ is a subset of the prime numbers.  We let
$\pi_\mathcal{P}(x)$ denote the number of primes $p\le x$ with $p\in\mathcal{P}$.
For a positive integer $n$ we let $n_\mathcal{P}$ denote the largest divisor
of $n$ that is free of prime factors outside of $\mathcal{P}$.

Let $e$ be an integer with $e>1$.  Let $\epsilon(x)$ be an arbitrary
monotonic function with
\begin{equation}
\label{epsilonspecs}
\epsilon(x)=o(1),~~~
\epsilon(x)>1/\log\log x,~~~ \epsilon(x^{1/\log\log x})<2\epsilon(x),
\end{equation}
where the last two conditions hold for $x$ sufficiently large.
We now partition the primes into 3 sets:
\begin{eqnarray*}
\LL&=&\{p~\text{prime}~:~\oo(e,p)\le p^{1/2}/\log p\}\\
\MM&=&\{p~\text{prime}~:~p^{1/2}/\log p<\ord(e,p)\le p^{1/2+2\epsilon(p)}\}\\
\HH&=&\{p~\text{prime}~:~\ord(e,p)>p^{1/2+2\epsilon(p)}\},
\end{eqnarray*}
where we use the mnemonic low, medium, high for $\LL, \MM, \HH$.
Note that $\LL$ contains the prime factors of $e$.

Let $\omega(n)$ denote the number of prime number divisors of $n$.

\begin{lem}
\label{lem:Lnumbers}
We have $\pi_\LL(x)=O(x/\log^3x)$ so that $\sum_{p\in\LL}1/p=O(1)$.
In addition, we have
\begin{equation}
\label{Lsum}
\sum_{n_\LL=n}\frac1n~=~\prod_{p\in\LL}(1-1/p)^{-1}~=~O(1)
\end{equation}
and
\begin{equation}
\label{Lcount}
\sum_{n_\LL=n,~n\le x}1~\ll~x/\log^3x.
\end{equation}
\end{lem}

\begin{proof}
To see the first assertion, let $y=x^{1/2}/\log x$ and note that if $p\in\LL$ and $p\le x$, 
then $\oo(e,p)\le y$.  That is, $p$ divides $e$ or some $e^j-1$ with
$1\le j\le y$.  Using the estimate $\omega(m)\ll\log m/\log\log m$, we have
\[
\pi_\LL(x)~\le~\omega\left(e\prod_{1\le j\le y}(e^j-1)\right)
~\ll~y^2/\log y~\ll~x/\log^3x.
\]
The result about $\sum_{p\in\LL}1/p$ then follows by partial summation,
and \eqref{Lsum} follows trivially as a consequence.  

We now prove \eqref{Lcount}.  Let $L_k(x)$ denote the number of integers
$n\le x$ with $n=n_\LL$ and $\omega(n)=k$.  We show by induction that
there is a positive constant $c$ such that
\begin{equation}
\label{hranalog}
L_k(x)~\le~c\frac{x}{(k-1)!\log^3x}\left(8\sum_{p\in\LL}\frac1{p-1}\right)^{k-1},
\end{equation}
 from which \eqref{Lcount} directly follows by summing on $k$ getting
\[
\sum_{n_\LL=n,~n\le x}1~\le~c\frac{x}{\log^3x}\exp\left(8\sum_{p\in\LL}\frac1{p-1}\right)
~\ll~\frac{x}{\log^3x}.
\]
To see \eqref{hranalog} note that we have already verified it in the case $k=1$.
Assume it is true at $k$.  Since no number can have two coprime prime-power divisors
bigger than the squareroot, we have
\begin{eqnarray*}
L_{k+1}(x)&\le&\frac1k\sum_{p\in\LL,~p^a\le x^{1/2}}L_k(x/p^a)\\
&\le&c\frac1{k!}\left(8\sum_{p\in\LL}\frac1{p-1}\right)^{k-1}
\sum_{p\in\LL,~p^a\le x^{1/2}} \frac{x/p^a}{\log^3(x/p^a)}\\
&\le&c\frac1{k!}\left(8\sum_{p\in\LL}\frac1{p-1}\right)^{k}\frac{x}{\log^3x}.
\end{eqnarray*}
This completes the proof of the lemma.
\end{proof}
Note that \eqref{Lsum} is all we shall need in this section, but we need the
stronger result \eqref{Lcount} for our later results.

For a positive integer $n$, let $\gamma(n)$ denote the largest squarefree
divisor of $n$, sometimes called the ``core" of $n$.  

\begin{lem}
\label{lem:trivpart}
But for a set of natural numbers $n$ of asymptotic density $0$ we have
\begin{eqnarray*}
n_\LL&<&\log n\\
n/\gamma(n)&<&\log n\\
\omega(n)&<&2\log\log n.
\end{eqnarray*}
\end{lem}

\begin{proof}
The first assertion follows directly from  \eqref{Lsum}.
The assertion about $n/\gamma(n)$ follows from the 
fact that the number of $n\le x$ with $n/\gamma(n)>T$ is
$O(x/\sqrt{T})$.  Indeed, if $u=n/\gamma(n)$, then $u\gamma(u)|n$ and
$u\gamma(u)$ is squareful (divisible by the square of each of its prime
factors).  The assertion then follows from partial summation and the
fact that the number of squareful numbers up to $x$ is $O(\sqrt{x})$.
The final assertion about $\omega(n)$ follows from the theorem of
Hardy and Ramanujan that the normal number of prime factors of $n$ is
$\log\log n$.
\end{proof}

One question of interest is how large can we expect $n_\MM$ to be for most
numbers $n$.  Since most numbers do not have a divisor
very near their square root, there is hope that this ingredient can be
used.  Erd\H os and Murty used this idea to show that
$\pi_\MM(x)=o(\pi(x))$ and Pappalardi and Indlekofer--Katai got more
quantitative versions of this result.  We state a consequence from the latter paper.
\begin{lem}[\cite{indlekofer-timofeev-divisors-of-shifted-primes},
 Cor. 6]
\label{lem:indlekofer-primes}
With $\epsilon(x)$ as specified in \eqref{epsilonspecs}, we have
$\pi_\MM(x)=O(\epsilon(x)^{1/12}\pi(x))$.
\end{lem}

We now show that as a consequence of Lemma \ref{lem:indlekofer-primes}
not many integers $n$ have a large divisor composed of primes from $\MM$.
Let $\Lambda$ denote the von Mangoldt function.

\begin{lem}
\label{lem:Mnumbers}
With $\epsilon(x)$ as specified in \eqref{epsilonspecs},
the number of integers $n\le x$ with $n_\MM>n^{1/3}$ is $O(\epsilon(x)^{1/12}x)$.
\end{lem}
\begin{proof}
We have
\[
\sum_{n\le x}\log n_\MM=
\sum_{n\le x}\sum_{\substack{d|n\\ d_\MM=d}}\Lambda(d)
~=~\sum_{\substack{d_\MM=d\\ d\le x}}\Lambda(d)\left\lfloor\frac{x}{d}\right\rfloor
~\le~x\sum_{\substack{p\in\MM\\ p\le x}}\frac{\log p}{p}+O(x).
\]
Now, using Lemma \ref{lem:indlekofer-primes} and \eqref{epsilonspecs},
\begin{eqnarray*}
\sum_{p\in\MM,~p\le x}\frac{\log p}{p}
&=&\frac{\log x}{x}\pi_\MM(x)+\int_2^x\frac{\log t-1}{t^2}\pi_\MM(t)\,dt\\
&\ll&\int_2^x\frac{\epsilon(t)^{1/12}}t\,dt+o(1)\\
&=&\int_2^{x^{1/\log\log x}}\frac{\epsilon(t)^{1/12}}t\,dt
+\int_{x^{1/\log\log x}}^x\frac{\epsilon(t)^{1/12}}t\,dt+o(1)\\
&\ll&\frac{\log x}{\log\log x}+\epsilon(x)^{1/12}\log x
~\ll~\epsilon(x)^{1/12}\log x.
\end{eqnarray*}
Thus,
\[
\sum_{n\le x}\log n_\MM~\ll~\epsilon(x)^{1/12}x\log x,
\]
so that the result follows readily.
\end{proof}

\begin{lem}
\label{lem:lambdaest}
For $x$ sufficiently large, the number of integers $n\le x$ with 
$\lambda(n)\le n\exp(-(\log\log n)^3)$ is at most $x/(\log x)^{10}$.
\end{lem}

\noindent
This result follows from Theorem 5 of \cite{pomerance-power-generator}.

We are now ready to prove the first part of Theorem \ref{thm:big-order-mod-n}.

\begin{thm}
\label{thm:ordn1}
Suppose $\epsilon(n)$ satisfies \eqref{epsilonspecs}.  But for a
set of integers $n$ of asymptotic density~$0$ we have
\[
\oo(e,n)~>~n^{1/2+\epsilon(n)}.
\]
\end{thm}

\begin{proof}
By Lemma \ref{lem:lambdaest} we may assume that $\lambda(n)>n\exp(-(\log\log n)^3)$.
Thus, from Lemma \ref{lem:ordineq} and Lemma \ref{lem:trivpart} we have
\begin{eqnarray*}
&&\oo(e,n)~>~\exp(-(\log\log n)^3)\prod_{p|n/n_\LL}\ord(e,n)\\
&&~~\ge~
\exp(-(\log\log n)^3)\prod_{p|n_\MM}(p^{1/2}/\log p)~\prod_{p|n_\HH}p^{1/2+2\epsilon(p)}\\
&&~~\ge
~\exp(-(\log\log n)^3-\omega(n)\log\log n)\gamma(n_\MM)^{1/2}\gamma(n_\HH)^{1/2+2\epsilon(n)}\\
&&~~\ge~\exp(-2(\log\log n)^3)n^{1/2}n_\HH^{2\epsilon(n)}.
\end{eqnarray*}
By Lemmas \ref{lem:trivpart} and \ref{lem:Mnumbers} we may also assume that
$n_\HH>n^{3/5}$.  Thus, our result follows from \eqref{epsilonspecs}.
\end{proof}

\section{The $1/2+\epsilon$ results}

We now consider analogs of Theorem \ref{thm:ordn1} in certain interesting
cases.  Say an infinite subset $\Ss$ of the natural numbers has property P 
``almost always" if 
\[
\sum_{\substack{s\in\Ss,~s\le x\\ s\text{ has property P}}}1
~\sim~\sum_{s\in\Ss,~s\le x}1~~\text{ as}~~x\to\infty.
\]
In this section P will be the property that $\oo(e,\lambda(n))>n^{1/2+\epsilon(n)}$.
That is, for $\epsilon(x)$ satisfying \eqref{epsilonspecs},
\[
\text{$n$ has property P$_\epsilon$: }~\oo(e,\lambda(n))>n^{1/2+\epsilon(n)}.
\]
Our goal of this section is to prove the following theorem, which comprises
the union of the first items of Theorems~\ref{thm:big-order-mod-p-1},
\ref{thm:main-epsilon}, and \ref{thm:lambda}.

\begin{thm}
\label{thm:Pepsilon}
If $\epsilon(x)$ satisfies \eqref{epsilonspecs} then the following sets
have property {\rm P}$_\epsilon$ almost always: the set of prime numbers, the set of
integers $n=pl$ where $p,l$ are primes with $p<l<2p$, and the set of all
natural numbers.
\end{thm}

We will need the following form of the Brun--Titchmarsh inequality (see
\cite{halberstam-richert}, Theorem 3.8):
\begin{lem} Suppose $k,l$ are coprime integers with $k>0$ and
let $\pi(x,k,l)$ be the number of primes $p \leq x$ such
that $p \equiv l~\text{\rm (mod }k\rm{)}$.  Then
$\pi(x,k,l) \ll \displaystyle{\frac{x}{\phi(k) \log(x/k)}}$
uniformly for $x>k$.
\end{lem}

We begin with an analog of Lemma \ref{lem:trivpart} for shifted primes.
\begin{lem}
\label{lem:p-1trivpart}
But for a set of prime numbers $p$ of relative density $0$ within the
set of all primes, we have
\begin{eqnarray*}
(p-1)_\LL&<&\log p\\
(p-1)/\gamma(p-1)&<&\log p\\
\omega(p-1)&<&2\log\log p
\end{eqnarray*}
\end{lem}
\begin{proof}
Using \eqref{Lcount} we have that
\[
\sum_{n=n_\LL,~n>T}\frac1n~\ll~\frac1{\log^2T}.
\]
Thus, by a trivial argument we may assume that $(p-1)_\LL<p^{1/2}$.
The Brun--Titchmarsh inequality and 
\eqref{Lcount} allow one to handle the remaining
cases where $(p-1)_\LL$ is between $\log p$ and $p^{1/2}$ as follows.
It suffices to show that 
\[
\sum_{n\ge\frac12\log x,~n=n_\LL}\pi(x,n,1)=o(\pi(x)),
\]
but the sum is $\ll\pi(x)\sum_{n\ge\frac12\log x,~n=n_\LL}1/\phi(n)$.
Using the well-known estimate $1/\phi(n)\ll(\log\log n)/n$,
we have our result from \eqref{Lcount}.
The argument for $(p-1)/\gamma(p-1)$ is similar, namely that a
trivial argument is used when $(p-1)/\gamma(p-1)$ is large and the
Brun--Titchmarsh inequality when it is small.  The final assertion
follows from the main result of \cite{erdos35} that the normal
number of prime factors of $p-1$ is $\log\log p$.
\end{proof}

We now turn our attention to an analog of Lemma \ref{lem:Mnumbers}
for shifted primes.  
\begin{lem}
\label{lem:p-1Mnumbers}
With $\epsilon(x)$ as specified in \eqref{epsilonspecs}, the number
of primes $p\le x$ with $(p-1)_\MM>p^{1/3}$ is $O(\epsilon(x)^{1/24}\pi(x))$.
\end{lem}
\begin{proof}  Using Brun's or Selberg's sieve (see \cite{halberstam-richert},
Theorem~2.4 or Theorem~3.12) we have that the
number of primes $p\le x$ with $p-1$ divisible by a prime $q>x^{1-\epsilon(x)^{1/24}}$
is
\[
\le~\sum_{a\le x^{\epsilon(x)^{1/24}}}\sum_{\substack{q\le x/a\\ aq+1\text{ prime}}}1
~\ll~\frac{x}{\log^2x}\sum_{a\le x^{\epsilon(x)^{1/24}}}\frac1{\phi(a)}
~\ll~\epsilon(x)^{1/24}\pi(x),
\]
where we have used the well-known result that $\sum_{a\le T}1/\phi(a)\sim c\log T$
for an appropriate constant $c$.
Thus, we may assume that $p-1$ has no prime factor larger than $x^{1-\epsilon(x)^{1/24}}$.
Trivially we may also assume that $p-1$ has no prime-power factor this large as well.
Letting $\sum'$ denoting a sum over primes with these conditions, we have
\begin{eqnarray*}
\sum_{p\le x}{}^{'}\log(p-1)_\MM
&=&\sum_{p\le x}{}^{'}\sum_{\substack{d|p-1\\ d_\MM=d}}\Lambda(d)\\
&=&\sum_{\substack{d_\MM=d\\ d\le x^{1-\epsilon(x)^{1/24}}}}\Lambda(d)\pi(x,d,1)\\
&\ll&\sum_{\substack{d_\MM=d\\ d\le x^{1-\epsilon(x)^{1/24}}}}\Lambda(d)
\frac{x}{\phi(d)\log(x/d)}\\
&\le&\sum_{\substack{d_\MM=d\\ d\le x^{1-\epsilon(x)^{1/24}}}}
\Lambda(d)\frac{x}{d\epsilon(x)^{1/24}\log x},
\end{eqnarray*}
the penultimate estimate coming from the Brun--Titchmarsh inequality.
Using the first two displays in the proof of Lemma~\ref{lem:Mnumbers}, we have
\[
\sum_{\substack{d_\MM=d\\ d\le x}}\frac{\Lambda(d)}{d}
~\ll~\epsilon(x)^{1/12}\log x,
\]
so that with the above estimate, we get that
\[
\sum_{p\le x}{}^{'}\log(p-1)_\MM
~\ll~\epsilon(x)^{1/24}x.
\]
The lemma follows readily.
\end{proof}

The proof of Theorem \ref{thm:Pepsilon} for the set of prime numbers now
follows directly from the proof of Theorem~\ref{thm:ordn1} where we
replace Lemmas~\ref{lem:trivpart} and \ref{lem:Mnumbers} with
Lemmas~\ref{lem:p-1trivpart} and \ref{lem:p-1Mnumbers}, respectively.
Note that we may continue to use Lemma~\ref{lem:lambdaest} since the
estimate for the exceptional set in that lemma is $o(\pi(x))$.

We next turn our attention to the set of numbers $pl$ where $p,l$ are
primes with $p<l<2p$.  Proving Theorem~\ref{thm:Pepsilon} for this set
is equivalent to showing that
\begin{equation}
\label{plepsilon}
\oo(e,\lambda(pl))~>~Q^{1+\epsilon(Q)}
\end{equation}
for all but $o(\pi(Q)^2)$ pairs of primes $p,l\le Q$.

We have from \cite{pomerance-power-generator}, Theorem 6, the following
result in analogy to Lemma~\ref{lem:lambdaest}:  
But for $o(\pi(Q)^2)$ pairs of primes $p,l\le Q$ we have 
\begin{equation}
\label{lamlampl}
\lambda(\lambda(pl))~>~pl/\exp(2(\log\log Q)^3).
\end{equation}

Note that
\begin{equation}
\label{ordprod}
\oo(e,[a,b])~\ge~\oo(e,a)\oo(e,b)\frac{\lambda([a,b])}{\lambda(a)\lambda(b)}.
\end{equation}
Indeed, letting $A=\oo(e,a), B=\oo(e,b)$ we have
\[
\oo(e,[a,b])~=~[A,B]~=~
\frac{AB}{(A,B)}~\ge~\frac{AB}{(\lambda(a),\lambda(b))},
\]
so that using $\lambda([a,b])=[\lambda(a),\lambda(b)]$, \eqref{ordprod} follows.
We Apply \eqref{ordprod} with $a=p-1,b=l-1$, where $p,l$ are
distinct primes.  As $\lambda([p-1,l-1])=\lambda(\lambda(pl))$, we get
\begin{equation}
\label{ordlampl}
\oo(e,\lambda(pl))~>~\oo(e,p-1)\oo(e,l-1)\frac{\lambda(\lambda(pl))}{pl}.
\end{equation}
So, to show \eqref{plepsilon}, we assume that \eqref{lamlampl} holds and
we apply \eqref{ordlampl}.   The result follows from the fact that the
set of primes has property P$_\epsilon$ almost always.  (To be perfectly
precise, we use that the set of primes has property P$_{2\epsilon}$ almost always.)

The third class of numbers in Theorem \ref{thm:Pepsilon}, namely,
the set of all numbers $n$, is more difficult.  We begin with
a new result: 

\begin{thm} [Martin--Pomerance \cite{martinpom}]
\label{lem:lamlam}
As $n\to\infty$ through a certain set of integers of asymptotic density $1$, we have 
\[
\lambda(\lambda(n))=n\cdot\exp(-(1+o(1))(\log\log n)^2\log\log\log n)
\]
Thus, $\lambda(\lambda(n))>n/\exp((\log\log n)^3)$ almost always.
\end{thm}

We now give the analog result to Lemmas \ref{lem:trivpart} and \ref{lem:p-1trivpart}.
\begin{lem}
\label{lem:lamtrivpart}
We have
\begin{eqnarray*}
\lambda(n)_\LL&<&\exp((\log\log n)^2)\\
\lambda(n)/\gamma(\lambda(n))&<&\log n\\
\omega(\lambda(n))&<&(\log\log n)^2
\end{eqnarray*}
almost always.
\end{lem}
\begin{proof}
We have
\[
\sum_{n\le x}\log\lambda(n)_\LL
~\le~\sum_{n\le x}\sum_{\substack{p^a\|\lambda(n)\\ p\in\LL}}\log p^a
~\le~\sum_{\substack{p^a\le x\\ p\in\LL}}\log p^a\sum_{\substack{n\le x\\ p^a|\lambda(n)}}1.
\]
If a prime power $p^a$ divides $\lambda(n)$ it must be the case that either
$n$ is divisible by some prime $q\equiv1~(\text{mod }p^a)$ or $p^{a+1}|n$.
As
\[
\sum_{\substack{q\le x\\ q\text{ prime}\\q\,\equiv\,1\,(\text{mod } d)}}\frac1q
~=~\frac{\log\log x+O(\log d)}{\phi(d)}
\]
uniformly for all integers $d\ge 2$ (see \cite{pomamicable}, Theorem 1 and Remark 1,
or Norton \cite{norton}), we have
\[
\sum_{\substack{n\le x\\ p^a|\lambda(n)}}1
~\le~\frac{x}{p^{a+1}}
+\sum_{\substack{q\le x\\ q\text{ prime}\\q\,\equiv\,1\,(\text{mod } p^a)}}\frac{x}q
~=~\frac{x\log\log x}{\phi(p^a)}+O\left(\frac{x\log p^a}{p^a}\right).
\]
Hence
\begin{eqnarray*}
\sum_{n\le x}\log\lambda(n)_\LL
&\ll&x\log\log x\sum_{\substack{p^a\le x\\ p\in\LL}}\frac{\log p^a}{p^a}
+x\sum_{\substack{p^a\le x\\ p\in\LL}}\frac{(\log p^a)^2}{p^a}\\
&\ll&x\log\log x,
\end{eqnarray*}
the last inequality coming from the estimate for $\pi_\LL(x)$ in Lemma \ref{lem:Lnumbers}.
Thus we immediately get the first assertion in the lemma.

For the second assertion note that from (6) and (7) in
\cite{eps} we have
\[
\log(\lambda(n)/\gamma(\lambda(n)))~\ll~
\log\log x/\log\log\log x
\]
for all but $o(x)$ choices of $n\le x$.  Thus we have the second assertion.

The third assertion follows from the fact that the normal order of
$\omega(\lambda(n))$ is $\frac12(\log\log n)^2$, see \cite{ep}.
\end{proof}

Now we give the analog result to Lemmas \ref{lem:Mnumbers} and \ref{lem:p-1Mnumbers}.
\begin{lem}
\label{lem:lamMnumbers}
Let $\epsilon(x)$ satisfy \eqref{epsilonspecs}.  Almost all numbers $n$ have
the property that $\lambda(n)_\MM<n^{2/5}$.
\end{lem}
\begin{proof}
Let
\[
\MM'~=~\{p\text{ prime}:(p-1)_\MM>p^{1/3}\}.
\]
Lemma \ref{lem:p-1Mnumbers} tells us that $\pi_{\MM'}(x)\ll\epsilon(x)^{1/24}\pi(x)$.
We apply the proof of Lemma~\ref{lem:Mnumbers} with $\MM$ replaced by $\MM'$
and with $\epsilon(x)^{1/12}$ replaced by $\epsilon(x)^{1/24}$.  Thus, by
the final display of Lemma~\ref{lem:Mnumbers} we have that
\[
\sum_{n\le x}\log n_{\MM'}~\ll~\epsilon(x)^{1/24}x\log x.
\]
We thus get that $n_{\MM'}\le n^{1/12}$ almost always.  Assume that $n$ has this
property.  By Lemma~\ref{lem:trivpart}, we may also assume that $n/\gamma(n)<n^{1/90}$.
Thus, 
\begin{eqnarray*}
\lambda(n)_\MM&\le&
(n/\gamma(n))\lambda(\gamma(n))_\MM
~<~n^{1/90}\prod_{p|n}(p-1)_\MM\\
&=&n^{1/90}\prod_{p|n_{\MM'}}(p-1)_\MM\prod_{p|n/n_{\MM'}}(p-1)_\MM\\
&\le&n^{1/90}\gamma(n_{\MM'})\gamma(n/n_{\MM'})^{1/3}
~\le~n^{1/90}\,n_{\MM'}^{2/3}\,n^{1/3}~\le~n^{2/5}.
\end{eqnarray*}
This completes the proof of the lemma.
\end{proof}

We are in a position now to complete the proof of Theorem \ref{thm:Pepsilon}.
Assume that $n$ satisfies
the properties in Theorem~\ref{lem:lamlam} and 
Lemmas~\ref{lem:lamtrivpart},~\ref{lem:lamMnumbers}.
By Lemma~\ref{lem:lambdaest} we may also assume that $\lambda(n)>n\exp(-(\log\log n)^3)$.
Thus, $\lambda(n)_\HH>n^{3/5}/\exp(2(\log\log n)^3)$.
Using Lemma~\ref{lem:ordineq} and assuming that $n$ is large, we have
\begin{eqnarray*}
\oo(e,\lambda(n))&\ge&\frac{\lambda(\lambda(n))}{\lambda(n)}\prod_{p|\lambda(n)}\oo(e,p)\\
&>&\exp(-(\log\log n)^3)\prod_{p|\lambda(n)_\MM}(p^{1/2}/\log p)
\prod_{p|\lambda(n)_\HH}p^{1/2+2\epsilon(p)}\\
&>&\exp(-2(\log\log n)^3)\gamma(\lambda(n)_\MM)^{1/2}\gamma(\lambda(n)_\HH)^{1/2+2\epsilon(n)}\\
&>&\exp(-3(\log\log n)^3)\lambda(n)^{1/2}\lambda(n)_\HH^{2\epsilon(n)}\\
&>&\exp(-4(\log\log n)^3)n^{1/2+(6/5)\epsilon(n)}\\
&>&n^{1/2+\epsilon(n)}.
\end{eqnarray*}
This completes the proof of Theorem \ref{thm:Pepsilon}.

\section{The $1/2+c$ results}

The spirit of Theorems \ref{thm:ordn1} and \ref{thm:Pepsilon}
concerns the best that can be said for almost all cases.
In this section we relax the ``almost all" to ``a positive proportion"
and so prove somewhat stronger results.  One could relax further
to ``infinitely often," but then it occurs that quite cheap results
can be had.  For example, if $p$ is a prime that does not divide $e$,
then $\ord(e,p^j)=p^{j-O(1)}$, so that $\ord(e,n)\gg n$ infinitely often.

We begin with the case of $\ord(e,p)$ for $p$ prime.  As mentioned in
the Introduction, one way of getting a fairly decent result here is to
have a very large prime factor of $p-1$ as afforded by a series of 
papers culminating in the recent paper \cite{baker-harman-shifted-primes}.

\begin{lem}[Baker--Harman]
\label{lem:bh}
For a positive proportion of the primes $p$, there is a prime $q|p-1$ with
$q>p^{0.677}$.
\end{lem}
\noindent
Note that this result follows from (7.1) in \cite{baker-harman-shifted-primes}.

We use this result to immediately get the following:
\begin{lem}
\label{lem:ordp}
We have $\ord(e,p)>p^{0.677}$
for a positive proportion of the primes $p$. 
\end{lem}
\begin{proof}
Among the primes $p$ for which $p-1$ is divisible by a prime $q>p^{0.677}$,
consider those for which $\ord(e,p)$ is not divisible by $q$.  Then if $p\le x$, we
have $\ord(e,p)<x^{0.323}$.  As in the argument
for $\pi_\LL(x)$ in the proof of Lemma~\ref{lem:trivpart}, the number
of such primes is $O(x^{0.646}/\log x)=o(\pi(x))$.  Thus, only a negligible
number of primes which satisfy the previous lemma do not satisfy the present lemma.
\end{proof}

Our basic strategy in this section to make
$\oo(e,m)$ large, is to manage to place in $m$ a large prime $p$ for which $\ord(e,p)$
is large, and then use the ideas of the previous sections to show that the remainder
of $m$ cannot do too much damage most of the time.  For $\oo(e,n)$ the idea
is especially transparent.

\begin{thm}
\label{thm:ordn2}
We have $\oo(e,n)>n^{0.677}$
for a positive proportion of integers $n$.
\end{thm}
\begin{proof}
The only subtlety here is that we need to extend Lemma \ref{lem:bh}
slightly.  By the Brun--Titchmarsh inequality, the proportion of primes
$p$ with a prime factor $q$ of $p-1$ in the interval $[p^{0.677},p^{0.677+2\epsilon}]$
is $O(\epsilon)$.  So if $\epsilon$ is small enough compared to the positive
proportion produced in Lemma~\ref{lem:bh}, then there must be a positive proportion
left over with $q>p^{0.677+2\epsilon}$.  And, for all but a negligible proportion
of these numbers, as in Lemma~\ref{lem:ordp}, we have $\ord(e,p)>p^{0.677+2\epsilon}$.
Now consider for such primes $p$, integers of the form $ap\le x$,
where $a\le x^{\epsilon}$.  For such primes $p\le x$ the number of integers $a$ that may
be taken is $\gg x/p$, and letting $p$ run from $x^{1-\epsilon}$ to $x$ there is
never any double counting of any $ap$.  Thus, the number of such numbers $ap$ is
$\gg\sum x/p\gg x$.  Further, 
\[
\oo(e,ap)\ge\ord(e,p)>p^{0.677+2\epsilon}>(ap)^{0.677}.
\]
This completes the proof of the theorem.
\end{proof}

We say $n$ has
property P$_c$ if $\oo(e,\lambda(n))>n^{1/2+c}$.  
In the rest of this section we take $c=0.092$.
\begin{thm}
\label{thm:Pc}
Positive proportions of the set of primes and the set of all natural numbers
have property {\rm P}$_c$.  Further, there are $\gg \pi(Q)^2$ pairs of primes
$p,l\le Q$ such that $pl$ has property {\rm P}$_c$.
\end{thm}
\begin{proof}
We begin with the case of primes, from which the other two cases will
follow easily.  We actually show a slightly stronger result: there is
some $\delta>0$ such that a positive proportion of the primes have
property P$_{c+\delta}$.  Let $\PP$ be the set of primes $q$ for which $\ord(e,q)>q^{0.677}$.
Lemma~\ref{lem:ordp} tells us that this set of primes comprises a positive proportion
of all primes.  Consider primes $p\le x$ where $q|p-1$ for
some $q\in\PP$ and with $x^{0.52-\epsilon}<q\le x^{0.52}$.  Here,
$\epsilon>0$ is arbitrarily small but fixed.  It follows from 
\cite{baker-harman-brun-titchmarsh-on-average}, Theorem 1, that a positive
proportion of primes $p$ are so representable.  Further, it follows from
Lemma~\ref{lem:p-1trivpart} that by neglecting only a relative density 0 of 
such primes $p$, we have 
\begin{eqnarray*}
\oo(e,p-1)&>&(p/q)^{1/2-o(1)}q^{0.677}~=~p^{1/2-o(1)}q^{0.177+o(1)}\\
&>&p^{1/2+(0.52-\epsilon)(0.177)-o(1)}.
\end{eqnarray*}
As $(0.52)(0.177)>c$, if $\epsilon$ is taken small enough, we have
$(0.52-\epsilon)(0.177)>c+\delta$ for some fixed $\delta>0$.  Thus,
$\oo(e,p-1)>p^{1/2+c+\delta}$, with this holding for a positive proportion of primes $p$.
Thus, we have the theorem for the set of primes.

Now consider the numbers $pl$, where $p,l$ are primes with $p,l\le Q$.
We apply \eqref{ordlampl} where $p,l$ are
primes with $p,l\le Q$ which have property P$_{c+\delta}$.  Assuming 
as we may that $pl$ satisfies \eqref{lamlampl}, we have
\[
\oo(e,\lambda(pl))~>~(pl)^{1/2+c+\delta}\exp(-2(\log\log Q)^3).
\]
Thus, there are $\gg\pi(Q)^2$ pairs of primes $p,l\le Q$ for which
$pl$ has property P$_c$.

We now consider the set of all
positive integers.  Consider the integers $n=ap$ where $a\le p^{\delta/2}$,
where $p$ is a prime with property P$_{c+\delta}$.  By the first part of the proof,
these numbers $n$ comprise a positive proportion of all numbers $n$.
Further, for such a number $n$ we have
\[
\oo(e,\lambda(n))~\ge~\oo(e,p-1)~>~p^{1/2+c+\delta}>(ap)^{1/2+c}~=~n^{1/2+c}.
\]
Thus, $n$ has property P$_c$.  This completes the proof of the theorem.

\end{proof}

\section{The $1-\epsilon$ results}

In this section we improve the $1/2+\epsilon$ results to $1-\epsilon$,
but we assume the Generalized Riemann Hypothesis (GRH).
We begin with the following slight strengthening of Theorem 2 of
\cite{artincat}: 
\begin{thm}
\label{thm:order-on-grh}
Let $e\ge2$ be an integer.
If the GRH is true, then for $x,y$ with $1\le y\le\log x$,
\[
\left|\left\{ p \leq x : \ord(e,p) \leq \frac{p}{y}\right \}\right|
~\ll~  \displaystyle{\frac{\pi(x)}{y}}
+ 
\frac{x \log \log x }{\log^2 x},
\]
where the implied constant depends at most on the choice of $e$.
\end{thm}
\begin{proof}  Since the proof is rather similar to the proof of
the main theorem in \cite{hooley-artin} and the proof of Theorem~2
in \cite{artincat}, we only give a brief outline.
With $i_p = 
(p-1)/\ord(e,p)$, we see that $\ord(e,p) \leq p/y$ implies that $i_p\ge y/2$.

{\em First step:}  We first consider primes $p$ such that 
$i_p \in ((x\log x)^{1/2}, x)$.  
As in the first part of the proof of Lemma~\ref{lem:Lnumbers},  the
number of such primes is $O( x/ \log^2 x) $.

{\em Second step:} Consider primes $p$ such that $q | i_p$ for some prime $q
$ in the interval $ [ \frac{x^{1/2}}{\log^3 x}, (x \log x)^{1/2} ]$. We
may bound this by considering primes 
$p \leq x$ such that $p \equiv  1~({\rm mod}~ q)$ for
some prime $q \in [ \frac{x^{1/2}}{\log^3 x}, (x \log x)^{1/2} ]$.
The Brun--Titchmarsh inequality then gives
that the number of  such primes $p$ is at most
\begin{equation*}
\sum_{ q \in [ \frac{x^{1/2}}{\log^3 x}, (x \log x)^{1/2} ]  }
\frac{x}{\phi(q) \log(x/q)}
~\ll
~\frac{x}{\log x}
\sum_{ q \in [ \frac{x^{1/2}}{\log^3 x}, (x \log x)^{1/2} ]  }
\frac{1}{q}
~\ll
~\frac{x \log\log x}{\log^2 x}.
\end{equation*}

{\em Third step:} Now consider primes $p$ such that $q| i_p$ for some prime
$q$ in the interval $[ y, \frac{x^{1/2}}{\log^3 x}]$.
In this range the 
GRH gives useful bounds; by (28) in \cite{hooley-artin} or Corollary~6 and
Lemma~9 of \cite{artincat}, we have 
$$
|\{ p \leq x : q \mid i_p  \}|
~\ll 
~\frac{\pi(x)}{q \phi(q)} + O( x^{1/2} \log( x q^2) ).
$$
Summing over $q$, we find that the number of such $p$  is bounded by 
\begin{equation*}
\sum_{ q \in [ y, \frac{x^{1/2}}{\log^3 x}]  }
\left( 
  \frac{\pi(x)}{q^{2}} + O( x^{1/2} \log( x q^2) )
\right)
~\ll
~\frac{\pi(x)}{y}
+
\frac{x}{\log^2 x}.
\end{equation*}

{\em Fourth step:} For the remaining primes $p$, any prime divisor
$q|i_p$ is smaller than $y$. Hence $i_p$ must be divisible by
some integer $d$ in the interval $[y/2,y^2]$.  The analog of
(28) in \cite{hooley-artin} for not-necessarily-squarefree integers,
or more directly, Corollary~6 and
Lemma~9 of \cite{artincat}, gives
\begin{equation}
\label{chebotarev}
|\{ p \leq x : d \mid i_p  \}|
~\ll 
~\frac{\pi(x)}{d \phi(d)} + O( x^{1/2} \log( x d^2) ).
\end{equation}
Hence the total number of such $p$  is bounded by 
\begin{equation*}
\sum_{ d \in [y/2,y^2]  }
\left( 
  \frac{\pi(x)}{d \phi(d)} + O( x^{1/2} \log( x d^2) )
\right)
~\ll
~\frac{\pi(x)}{y},
\end{equation*}
where the last estimate follows from the well-known result $\sum_{a\le
  T}1/\phi(a)=c\log T+O(1)$ (for 
an appropriate constant $c$) and partial summation.
\end{proof}
\noindent{\em Remark.}  It follows easily from \eqref{chebotarev} that
for $1\le y\le x^{1/4}/\log x$ and assuming the GRH, we have 
\[
\left|\left\{ p \leq x : p^{1/2}y\log^2x\le \ord(e,p) \leq \frac{p}{y}\right \}\right|
~\ll~  \displaystyle{\frac{\pi(x)}{y}}.
\]

Let $\delta(x)=\sqrt{\log\log x/\log x}$.
By a slight abuse of notation, say an integer $n$ has property P$_{1-\delta}$
if $\oo(e,\lambda(n))\ge n^{1-\delta(n)}$.
Theorem \ref{thm:order-on-grh} is our principal tool in the proof of
the following result.  
\begin{thm}
\label{thm:grhthm}
Assume the GRH holds.
The set of primes and the set of integers $pl$ with $p,l$ prime and $p<l<2p$
have property {\rm P}$_{1-\delta}$ almost always.
\end{thm}
\begin{proof}
Let
\[
\WW~=~\{p\text{ prime}:\oo(e,p)<p/\log p\},
\]
where we use the mnemonic $\WW$ for weak.  From Theorem \ref{thm:order-on-grh}
we have
\begin{equation}
\label{piW}
\pi_\WW(x)~\ll~x\log\log x/\log^2 x.
\end{equation}
We now consider
\[
S:~=~\sum_{p\le x}\log(p-1)_\WW\,,
\]
following the lines of the proof of Lemma \ref{lem:p-1Mnumbers}.  We have
\begin{equation}
\label{p-1sum}
S~=~\sum_{\substack{d\le x\\ d_\WW=d}}\Lambda(d)\pi(x,d,1)
~=~\sum_{\substack{p\le x\\ p\in\WW}}\pi(x,p,1)\log p+O\left(\frac{x}{\log x}\right).
\end{equation}
Using Brun's or Selberg's sieve as in the proof of Lemma \ref{lem:p-1Mnumbers}, we have 
$\sum_{p>x^{1-\epsilon}}\pi(x,p,1)\ll\epsilon x/\log x$, so that
the contribution to the last sum in \eqref{p-1sum} from the primes $p>x^{1-\epsilon}$ is
$\ll \epsilon x$.  For primes $p\le x^{1-\epsilon}$ we use the Brun--Titchmarsh
inequality to get $\pi(x,p,1)\ll x/(\epsilon p\log x)$, so that using
\eqref{piW}, the contribution to the sum from these primes is $\ll x/(\epsilon\log x)$.
Letting $\epsilon=1/\sqrt{\log x}$, we get
\begin{equation}
\label{p-1Wnumbers}
\sum_{p\le x}\log(p-1)_\WW~\ll~x/\sqrt{\log x}.
\end{equation}
Thus, $(p-1)_\WW\le p^{\delta(p)/2}$ almost always.  The proof of our theorem
for the set of primes now follows in exactly the same way as in Theorem~\ref{thm:Pepsilon}.

The case for the numbers $pl$ now also follows using \eqref{lamlampl}
and our prior arguments.
\end{proof}

%
We now begin to examine the normal contribution to $\lambda(n)$ from primes in $\WW$.
\begin{lem} 
\label{l:no-big-p}
Assuming the GRH is true, for $x,\bigA\ge3$, the
number of integers $n \leq x$ such that $p|\lambda(n)$ for $p \in \WW$ and
$p>\bigA$ is 
$$
\ll~ x \log \log x \cdot \frac{\log \log \bigA }{\log \bigA}.
$$
\end{lem}
\begin{proof}
If $p|\lambda(n)$, then either $p^2 | n$ or some prime $q
\equiv 1~({\rm mod}~p)$ divides $n$.  The number of $n \leq x$ in the first
case is clearly bounded by $x/\bigA$.  
By the Brun--Titchmarsh
inequality and partial summation,
$$
x\sum_{q\le x,\,q\equiv1\,({\rm mod}\,p)}\frac1q
~\ll~\frac{x\log\log x}{p},
$$
hence the number of $n \leq x$ for which the second case occurs is
$$
\ll~\sum_{p > \bigA, \,p \in \WW}
\ 
\sum_{q \leq x, \,q \equiv 1\,({\rm mod}\,p)}
x/q
~\ll 
~x \log \log x \sum_{p > \bigA, \,p \in \WW} 1/p
$$
which, since $\pi_{\WW}(x) \ll x \log \log x/\log^2 x$, is 
$$
\ll
~x\log \log x \cdot \frac{\log \log \bigA}{\log \bigA}
$$
by partial summation.
\end{proof}
%
We now prove that for most integers $n$, $\lambda(n)_\WW$
is fairly small in the following sense:
\begin{lem}
\label{lem:fn-average}
Let $f(n)=\sum_{p|\lambda(n),\,p\in\WW}\log p$. 
Assuming the GRH is true, for almost all integers $n$, we have
$$
f(n)~<~ (\log \log n)^2  .
$$  
\end{lem}
\begin{proof}
Take
$
\bigA = 
\exp \left( (\log \log x)^2 \right)
$
in Lemma~~\ref{l:no-big-p}.  Then the number of $n \leq x $ for which 
some $p \in \WW$, $p> \bigA $ divides $\lambda( n)$ is $o(x)$.
Letting $\sum'$ denote a sum over
$n$ for which no $p \in 
\WW$, $p> \bigA $ divides $\lambda( n)$, we obtain as before that
\begin{eqnarray*}
\sum_{n \leq x}{}^{'} f(n)
&=&\sum_{p\le T,\,p\in\WW}\log p
\sum_{\substack{n\le x \\p|\lambda(n)}}\hspace{-1mm}{}^{'} \, 1 \\
&\ll&
x \sum_{p \leq T, \,p \in \WW} \frac{\log p}{p^2}
+ x \log \log x 
\sum_{p \leq \bigA, \,p \in \WW} \frac{\log p}{p}.
\end{eqnarray*}
Since $\pi_{\WW }(x) \ll x\log\log x/\log^2 x$, partial summation gives that
$$
\sum_{p \leq T, \,p \in \WW } \frac{\log p}{p} 
~\ll
~(\log \log T)^2~\ll~(\log\log\log x)^2.
$$
Hence
$$
\sum_{n \leq x}{}^{'} f(n)
~\ll ~x \log \log x\, (\log \log \log x)^2.
$$
Thus, the average order of $f(n)$, after removing those integers $n$ where
$\lambda(n)$ is divisible by some $p \in \WW$, $p>T$, 
is $\ll \log\log n(\log\log\log n)^2$.  We conclude that
$f(n)<(\log \log n)^2$ holds for almost all $n$.
\end{proof}

We are now ready to prove a result for $\oo(e,\lambda(n))$ on the
assumption of the GRH.  
\begin{thm}
\label{thm:stronglambda}
If the GRH is true, then for each fixed integer $e\ge2$,
\[
\oo(e,\lambda(n))~=~n\cdot\exp(-(1+o(1))(\log\log n)^2\log\log\log n)
\]
as $n\to\infty$ through a set of asymptotic density $1$.
\end{thm}

\begin{proof}
We shall show that if the GRH is true, then
\begin{equation}
\label{stronglower}
\oo(e,\lambda(n))\ge
\lambda(\lambda(n))\exp\left(-3(\log\log n)^2(\log\log\log\log n)^2\right)
\end{equation}
for almost all $n$.  The theorem will then follow from the trivial inequality
$\oo(e,\lambda(n))\le\lambda(\lambda(n))$ and Theorem \ref{lem:lamlam}.
By Lemma~\ref{lem:fn-average} we may assume that
$f(n)<(\log\log n)^2$. 
Let
\[
\WW_1~=~\{p\text{ prime}:p/\log p~\le ~\oo(e,p)
<~p/(\log\log p\cdot\log\log\log p)\},
\]
so that by Theorem \ref{thm:order-on-grh} we have 
$\pi_{\WW_1}(x)\ll\pi(x)/(\log\log x\cdot\log\log\log x)$.
Let $g(n)=\sum_{p|\lambda(n),\,p\in\WW_1}1$.  Then
\begin{eqnarray*}
\sum_{n\le x}g(n)&=&\sum_{p\le x,\,p\in\WW_1}\sum_{n\le x,\,p|\lambda(n)}1\\
&\ll&x\sum_{p\le x,\,p\in\WW_1}\frac1{p^2}
+x\log\log x\sum_{p\le x,\,p\in\WW_1}\frac1p\\
&\ll&x\log\log x\cdot\log\log\log\log x,
\end{eqnarray*}
the last estimate coming from partial summation and our inequality for $\pi_{\WW_1}(x)$.
Thus, for almost all $n$, $g(n)<\log\log n\,(\log\log\log\log n)^2$.

Also, let
\begin{eqnarray*}
\WW_2~=~\{p\text{ prime}:p/(\log\log p\cdot \log\log\log p )~\le ~\oo(e,p)&&\\
<~p/\log\log\log p\},&&
\end{eqnarray*}
so that by Theorem \ref{thm:order-on-grh} we have $\pi_{\WW_2}(x)\ll\pi(x)/\log\log\log x$.
We let $h(n)=\sum_{p|\lambda(n),\,p\in\WW_2}1$.  As in the calculation for $g(n)$,
we get
\[
\sum_{n\le x}h(n)~\ll~x(\log\log x)^2/\log\log\log x,
\]
so that for almost all $n$ we have 
\[
h(n)<(\log\log n)^2\log\log\log\log n/\log\log\log n.
\]

Now assume that $f(n),g(n),h(n)$ are bounded as above, and assume that the
inequalities in Lemma \ref{lem:lamtrivpart} hold.  We have by Lemma \ref{lem:ordineq}
\begin{equation}
\label{lamstart}
\oo(e,\lambda(n))~\ge~\frac{\lambda(\lambda(n))}{\lambda(n)}
\prod_{p|\lambda(n)}\oo(e,p)~\ge~\frac{\lambda(\lambda(n))}{\lambda(n)}ABC,
\end{equation}
where
\begin{eqnarray*}
A:&=&\prod_{p|\lambda(n)_{\WW_1}}\frac{p}{\log p}\,,\\
B:&=&\prod_{p|\lambda(n)_{\WW_2}}\frac{p}{\log\log p \cdot \log\log\log p}\,,\\
C:&=& \prod_{p|\lambda(n)/\lambda(n)_{\WW\cup\WW_1\cup\WW_2}}
\frac{p}{\log\log\log p}\,.
\end{eqnarray*}
Now
\[
ABC~\ge~\frac{\prod_{p|\lambda(n)/\lambda(n)_{\WW}}p}{DEF},
\]
where
\begin{eqnarray*}
D:&=&(\log n)^{g(n)},\\
E:&=&\left(\log\log n \cdot \log\log\log n \right)^{h(n)},\\
F:&=&(\log\log\log n)^{\omega(\lambda(n))}.
\end{eqnarray*}
By our assumptions on $n$, and taking $n$ sufficiently large, we have
\[
DEF~\le~\exp(2(\log\log n)^2(\log\log\log\log n)^2).
\]
Further,
\[
\prod_{p|\lambda(n)/\lambda(n)_{\WW}}p
~=~\frac{\gamma(\lambda(n))}{\exp(f(n))}
~\ge~\frac{\lambda(n)}{\log n \cdot\exp\left((\log\log n)^2\right)}.
\]
Hence by our above estimates,
\[
ABC
~\ge~\lambda(n)\exp\left(-3(\log\log n)^2(\log\log\log\log n)^2\right)
\]
for almost all $n$.  We use this estimate in \eqref{lamstart}, so that \eqref{stronglower}
and the theorem follow.
\end{proof}

As mentioned in the introduction, $\oo(e,\lambda(n))$ is the period
of the power generator $u^{e^i}~(\text{mod }n)$ if $\oo(u,n)=\lambda(n)$,
that is, if $\oo(u,n)$ is as large as possible.  We now briefly consider
the situation for a general modulus $n$ when we do not make this assumption
about $u$.  
We have the following result.
\begin{thm}
\label{thm:per}
Assuming the GRH, for any fixed integers $e,u\ge2$, 
the period of the
sequence $u^{e^i}~(\text{\rm mod }n)$ is equal to
\[
n\cdot\exp(-(1+o(1))(\log\log n)^2\log\log\log n)
\]
as $n\to\infty$ through a certain set of integers
of asymptotic density~$1$.
\end{thm}
\begin{proof}
First note the elementary inequality
\begin{equation}
\label{elemineq}
\text{for }j\mid n\text{ we have }\oo(e, n/j)\ge\frac1j\oo(e,n).
\end{equation}
To see this, as before let $j_{(e)},n_{(e)}$ be the largest divisors of
$j,n$ respectively that are coprime to $e$, so that $\oo(e,n)=\ord(e,n_{(e)})$ and
$\oo(e,n/j)=\ord(e,n_{(e)}/j_{(e)})$.  Let $j_{(e)}=j_1j_2$ where $j_1$ is the largest divisor
of $j_{(e)}$ that is coprime to $n_{(e)}/j_{(e)}$.  Then
\[
\ord(e,n_{(e)})~=~\ord(e,j_1j_2n_{(e)}/j_{(e)})~=~[\ord(e,j_1),\ord(e,j_2n_{(e)}/j_{(e)})].
\]
Further, $\ord(e,j_2n_{(e)}/j_{(e)})\mid j_2\cdot\ord(e,n_{(e)}/j_{(e)})$, so that
\begin{eqnarray*}
\oo(e,n)~=~\ord(e,n_{(e)})&\le&\ord(e,j_1)\cdot j_2\cdot\ord(e,n_{(e)}/j_{(e)})\\
&\le&j_{(e)}\cdot\ord(e,n_{(e)}/j_{(e)})~\le j\cdot\oo(e,n/j),
\end{eqnarray*}
which proves \eqref{elemineq}.  Recall that the period for the
sequence $u^{e^i}~(\text{mod }n)$ is $\oo(e,\oo(u,n))$.
Thus, if $\oo(u,n)=\lambda(n)/j$, we have by \eqref{elemineq} that the period is
\[
\oo(e,\lambda(n)/j)~\ge~\frac1j\oo(e,\lambda(n)).
\]
But, on the GRH we have $\oo(u,n)>n/(\log n)^{2\log\log\log n}$
almost always; this follows from the proof of Cor.~2 in
\cite{li-pomerance-artin-for-composite}.
Thus, we may take $j<(\log n)^{2\log\log\log n}$, so the result
follows from Theorem \ref{thm:stronglambda}.
\end{proof}

\bibliographystyle{abbrv} 

\bigskip
\begin{thebibliography}{10}

\bibitem{baker-harman-brun-titchmarsh-on-average}
R.~C. Baker and G.~Harman.
\newblock The {B}run-{T}itchmarsh theorem on average.
\newblock In {\em Analytic number theory, Vol. 1 (Allerton Park, IL, 1995)},
  volume 138 of {\em Progr. Math.}, pages 39--103. Birkh\"auser Boston, Boston,
  MA, 1996.

\bibitem{baker-harman-shifted-primes}
R.~C. Baker and G.~Harman.
\newblock Shifted primes without large prime factors.
\newblock {\em Acta Arith.}, 83(4):331--361, 1998.

\bibitem{erdos35}
P. Erd\H os,
\newblock On the normal number of prime factors of $p-1$ and some
other related problems concerning Euler's $\varphi$-function.
\newblock {\em Quart.\ J. Math.\ (Oxford Ser.)} 6: 205--213, 1935.

\bibitem{erdos-murty}
P.~Erd{\H{o}}s and M.~R. Murty.
\newblock On the order of $a\pmod p$.
\newblock In {\em Number theory (Ottawa, ON, 1996)}, pages 87--97. Amer. Math.
  Soc., Providence, RI, 1999.

\bibitem{ep}
P. Erd\H os and C. Pomerance.
\newblock On the normal number of prime factors of $\varphi(n)$.
\newblock {\em Rocky Mountain J. Math.}, 15: 343--352, 1985.

\bibitem{eps}
P. Erd\H os, C. Pomerance, and E. Schmutz.
\newblock Carmichael's lambda function.
\newblock {\em Acta Arith.}, 58: 363--385, 1991.

\bibitem{pomerance-power-generator}
J.~B. Friedlander, C.~Pomerance, and I.~E. Shparlinski.
\newblock Period of the power generator and small values of {C}armichael's
  function.
\newblock {\em Math. Comp.}, 70(236):1591--1605, 2001.
Corrigendum.  {\em Math. Comp.}, 71(240):1803--1806, 2002.

\bibitem{halberstam-richert}
H.~Halberstam and H.-E. Richert.
\newblock {\em Sieve methods}.
\newblock Academic Press [A subsidiary of Harcourt Brace Jovanovich,
  Publishers], London-New York, 1974.
\newblock London Mathematical Society Monographs, No. 4.


\bibitem{hooley-artin}
C.~Hooley.
\newblock On {A}rtin's conjecture.
\newblock {\em J. Reine Angew. Math.}, 225:209--220, 1967.

\bibitem{indlekofer-timofeev-divisors-of-shifted-primes}
K.-H. Indlekofer and N.~M. Timofeev.
\newblock Divisors of shifted primes.
\newblock {\em Publ. Math. Debrecen}, 60(3-4):307--345, 2002.

\bibitem{artincat}
P.~Kurlberg.
\newblock On the order of unimodular matrices modulo integers.
\newblock {\em Acta Arith.}, 110(2):141--151, 2003.

\bibitem{cat2}
P.~Kurlberg and Z.~Rudnick.
\newblock On quantum ergodicity for linear maps of the torus.
\newblock {\em Comm. Math. Phys.}, 222(1):201--227, 2001.

\bibitem{li-pomerance-artin-for-composite}
S.~Li and C.~Pomerance.
\newblock On generalizing {A}rtin's conjecture on primitive roots to composite
  moduli.
\newblock {\em J. Reine Angew. Math.}, 556:205--224, 2003.

\bibitem{martinpom}
G. Martin and C. Pomerance.
\newblock The iterated Carmichael $\lambda$ function
and the number of cycles of the power generator.
\newblock To appear.

\bibitem{norton}
K. Norton.
\newblock  On the number of restricted prime factors of an integer. I.
\newblock {\em Illinois J. Math.}, 20(4), 681--705, 1976.

\bibitem{pappalardi-p-1-div}
F.~Pappalardi.
\newblock On the order of finitely generated subgroups of ${\bf {Q}}\sp *\pmod
  p$ and divisors of $p-1$.
\newblock {\em J. Number Theory}, 57(2):207--222, 1996.

\bibitem{pomamicable}
C. Pomerance.
\newblock On the distribution of amicable numbers.
\newblock {\em J. Reine Angew.\ Math.}, 293/294: 217--222, 1977.

\bigskip
\bigskip

\end{thebibliography}

\end{document}